\documentclass{article}
\usepackage{graphicx}
\usepackage{latexsym}
\usepackage{amsmath}
\usepackage{amsfonts}
\usepackage{float}
\usepackage{subcaption}
\usepackage{mathtools}
\usepackage[utf8]{inputenc}
\usepackage[T1]{fontenc}
\usepackage{euscript,fullpage,amsmath,amssymb}
\usepackage[T1]{fontenc}
\usepackage{enumerate}
\usepackage{graphicx}
\usepackage{float}
\usepackage{setspace}
\usepackage{amssymb}
\usepackage{bbm}
\usepackage{subcaption}
\usepackage{stmaryrd}
\usepackage{amsmath}
\usepackage[T2A,T1]{fontenc}
\usepackage[russian,english]{babel}

\newtheorem{proposition}{Proposition}
\newenvironment{proof}[1][Proof]{\noindent\textbf{#1.} }{\ \rule{0.5em}{0.5em}}
\newtheorem{example}{Example}
\begin{document}
\title{Matching of observations}
\author{Th\'eophile Caby,\\
 LAMIA, Université des Antilles, Fouillole, Guadeloupe\\
 caby.theo@gmail.com}
\date{}

\maketitle
\begin{abstract}
We study the statistical distribution of the closest encounter between observations computed along different trajectories of a mixing dynamical system. At the limit of large trajectories, the distribution is of Gumbel type and depends on the length of the trajectories and on the Generalized Dimensions of the image measure. It is also modulated by an Extremal Index, for which we give a formula in the case of expanding maps of the interval and regular observations. We discuss the implications of these results for the study of physical systems.
\end{abstract}

\section{Introduction}
Recently, the problem of the shortest distance between orbits of a dynamical system has gained interest. Barros, Liao and Rousseau give the asymptotic behavior of the shortest distance between two orbits, and show its connection with sequence matching problems \cite{short}. In a later publication, these results were generalized to multiple orbits and observed orbits \cite{encoded}. Meanwhile, in a series of paper, the distributions the times of first synchronization and of the closest encounter between several trajectories were also studied, using tools of Extreme Value Theory \cite{synchro,d2,dq}. This approach provides numerical methods to compute fractal dimensions and hyperbolicity indices associated with the system, but necessitate to work with real trajectories, while often, physicists only have access to observations collected along the trajectories of the studied system. In the present paper, we generalize this approach to such observed trajectories, and show that for observations with large enough rank, one can recover information on the underlying system.

\section{The general approach}
Let us consider the dynamical system $(M,T,\mu)$. We take $T:M \to M$ to be a discrete transformation (it could be a discretized version of a flow) that leaves the probability measure $\mu$ invariant. To model the process of measurement, we consider the $C^1$ function $f:M \to J$ that we call the {\em observation}. Both the phase space $M$ and the observation space $J$ are compact metric spaces endowed with two distances that we will both call $d$ to simplify notations. We generally take $J\subset\mathbb{R}^m$, as observational data are usually a collection of real numbers that can be arranged into vectors. Because we are interested in the statistical properties of observations, we need a measure that is supported in the observational space; the image measure $\mu_f$, defined by 
$$\mu_{f}(A)=\mu(f^{-1}(A)),$$ for all $A \subset J$ such that $f^{-1}(A)$ is $\mu-$measurable.\\

For our purpose, we define the following process:

$$
Y_i=-\log(\max_{j=2,\dots,q}d(f(T^ix_1),f(T^ix_j))),
$$

being $x_1,...,x_q \in M$ $q$ starting points drawn independently from the invariant probability measure $\mu$.\\

Let $s\in \mathbb{R}$. To follow the usual procedure in Extreme Value Theory, we consider a sequence of thresholds $u_n(s)$ such that 

\begin{equation}\label{tau}
\mu_q(Y_0 > u_n(s)) \sim \frac{e^{-s}}{n},
\end{equation}

where $\mu_q$ is the product measure with support in $M^q$.\\
 
We notice that, since the $q$ trajectories are independent, we also have:

\begin{equation}\label{dqf}
\begin{aligned}
\mu_q(Y_0> u_n(s)) &=\int_J  \mu_f (B(y,e^{-u_n}))^{q-1}d \mu_f(y)\\
                                                       & \sim e^{-u_nD^f_q(q-1)},
\end{aligned}
\end{equation}

being $B(y,r)$ a ball centered at $y\in J$ of radius $r$ and $D^f_q$ the generalized dimension of order $q$ of the image measure, defined for $q\neq 1$ by

\begin{equation}
D_q^f=\underset{r\to 0}{\lim} \frac{\log \int_{J} f_*\mu(B(x,r))^{q-1}df_*\mu(x)}{(q-1)\log r}.
\end{equation}

We will place ourselves in physical situations, where the limit defining the previous quantity exists.\\

To satisfy both scalings \ref{tau} and \ref{dqf}, we take 

$$u_n(s)=\frac{\log n}{D_q^f(q-1)}+\frac{s}{D_q^f(q-1)}.$$

We now consider the variable

$$M_n(x_1,...,x_q)=\max \{Y_0,\dots,Y_{n-1}\}.$$

The distribution of this maximum gives us the hitting time statistics in the target set 

$$S^q_n=\{(s_1,...,s_q)\in M^q, \max_{j=2,\dots,q}d(f(s_1),f(s_j)) < e^{-u_n}\},$$

that is when all the observations lay in the same ball of radius $e^{-u_n}$ for the first time. We now apply techniques from Extreme Value Theory, in particular the spectral theory of Keller and Liverani \cite{kl,k} which yields the asymptotic distribution of the maximum for a large class of exponentially-mixing systems and for {\em regular} observations. The cumulative distribution of the maximum

$${F_n}(u_n) = \mu_q(\{(x_1,...,x_q) \in M^q \mbox{ s.t. } M_n(x_1,...,x_q) \leq u_n  \})$$

converges, in the sense that 

\begin{equation}
|F_n(u_n(s)) -\exp(-\theta^f_q e^{-s})| \underset{n\to\infty}\to 0.
\end{equation}
The term $\theta^f_q$ is called the extremal index, and is a number comprised between 0 and 1 that quantifies the tendency of the process $(Y_i)$ to form clusters of high values. The spectral theory demands that the system is rapidly mixing and that the measure of the target sets $S^q_n$ goes to zero in a regular fashion. More detailed presentations of the theory and its domain of application can be found in various publications \cite{kl,k,synchro,d2,dq}. The theory is proven to be particularly applicable to expanding maps of the interval \cite{synchro} and certain well-behaved 2-dimensional systems \cite{bakersandro}. It is not our goal to give conditions of existence of the extreme value law that are more adapted to the present case, since these are in practice difficult to check in dimension more than 1 or 2. We will however provide numerical evidence of the convergence to the extreme value law. We will now discuss on the values of the different parameters of the limit law, that can acquire physical meaning. 

\section{The Generalized Dimensions of the image measure $D_q^f$}

We have seen in preceding section that the $D_q^f$ spectrum modulates the synchronization properties of the observations. In fact, these quantities play a central role in many of the statistical properties of the system. We can show that if the $D_q^f$ spectrum is well defined, it also influences the behavior of diverse local dynamical quantities associated with the observations. It is well known that both return and hitting times of a chaotic system in small balls (in fact a re-scaled version of these quantities) have large deviations that are governed by the spectrum of generalized dimensions of the invariant measure \cite{dq,ldr}. A similar relation also holds for (finite-resolution) local dimensions (see \cite{ldr,dq} and \cite{these} for discussion). This kind of large deviations relations are known to hold for real trajectories, but they also apply to the recurrence times of observations (and to the local dimension of the image measure). The rate function is now entirely determined by the spectrum of generalized dimensions of the image measure. To see this, one can perform a direct adaptation of the proofs in \cite{ldr} and \cite{dq} to the case of observations. In particular, hypothesis {\bf A-1} in \cite{dq} is satisfied for exponentially-mixing systems also when considering observed trajectories \cite{obsrec,jerobs}. Hypothesis {\bf A-2} is satisfied for ergodic systems for which the local dimensions are well defined and for typical observations. Indeed, in that case, Young's theorem \cite{young} insures the exact dimensionality of the underlying system, then theorem 4.1 in \cite{hk} insures the exact dimensionality of the image measure. Hypothesis {\bf A-4} concerns the existence and analyticity of the generalized dimensions of the image measure. We now investigate this matter.\\

In \cite{hk}, Hunt and Kaloshin give results concerning the effect of projection on the generalized dimensions for $1\ge q \ge 2$. In this range, they show that if $M$ is a compact subset of $\mathbb{R}^n$ and $J=\mathbb{R}^m$, and if the generalized dimension of order $q$, $D_q$($=D_q^{Id}$) of the invariant measure exists, then

\begin{equation}\label{hk}
D_q^f=\min(D_q,m)
\end{equation}

 for a prevalent set of $C^1$ observables. See \cite{prev} for a review of prevalence, which is a notion of genericity for infinite dimensional spaces. Very little is known for the behavior of $D_q^f$ outside this range, although still in \cite{hk}, the authors show that no analogous general result hold when $q >2$.\\
 
A first application of the theory is the computation of the correlation dimension ($D_2$) of a physical system using observations, which can be performed by fitting the empirical distribution of $M_n$ and extracting the desired parameter (this type of methods of computation of fractal dimensions has been widely used in climate lately \cite{nature,messori,d2,dq}). Indeed, from result \ref{hk}, $D_2^f=D_2$ whenever the observable $f$ is typical and has a large enough rank (larger than $D_2$). The latter can be obtained by recording simultaneously the system at different locations in space (by using gridded observables) or by considering delay coordinates observables used in embedding techniques \cite{takens}. Notice that it is here enough to take $m\ge D_2$ delay coordinates to access the correlation dimension $D_2$, and not the $\lceil2D_0\rceil$ required to reconstruct the attractor \cite{takens}.\\

  To investigate on the value of $D_q^f$ for $q>2$, we compare in figure \ref{dqf} the numerical estimates of $D_q^f$ for different observations $f$ and the generalized dimensions associated with a motion on a Sierpinski gasket, for which explicit formulas are known (see \cite{dq} for the presentation of the map and explicit formulae). The estimates are obtained by fitting the empirical distribution of the maximum value taken by the process $(Y_i)$ over blocks of size $5.10^4$. This procedure will also allow us to confirm the convergence of the distribution. The results are averaged over 10 runs, using trajectories of length $2.10^8$. The error bars represent the standard deviations of the results. Functions $f_1$, $f_2$ are diffeomorphisms, which are known to preserve the generalized dimensions. Indeed, for these two functions, good agreement is found, so that the two curves are hardly distinguishable visually in the figure. These results suggest that this method of computation of $D_q$ can be completed and even improved by introducing a diffeomorphism computed along the orbit of the system, which may, if well chosen, speed up the convergence of the method and provide better estimates. Function $f_3$ is a very oscillatory function, which gives a point in the observational space many antecedent, having the effect to alter significantly the fine structure of the image measure. We do not know whether the disagreement with the $D_q$ spectrum is due to the method not being at convergence, or if is a sign that the spectrum is not preserved under the action of $f_3$. However, the small disagreement for $q=2$ seems to imply that the method may not be at convergence, since the correlation dimension is preserved by typical observations. $f_4$ is not a diffeomorphism either, but has a more simple structure. For this function, the generalized dimensions seem to be preserved. $f_5$ is a degenerate function yielding values close to 1.\\
 \begin{figure}[h!]
 \centering
\includegraphics[height=3.5in]{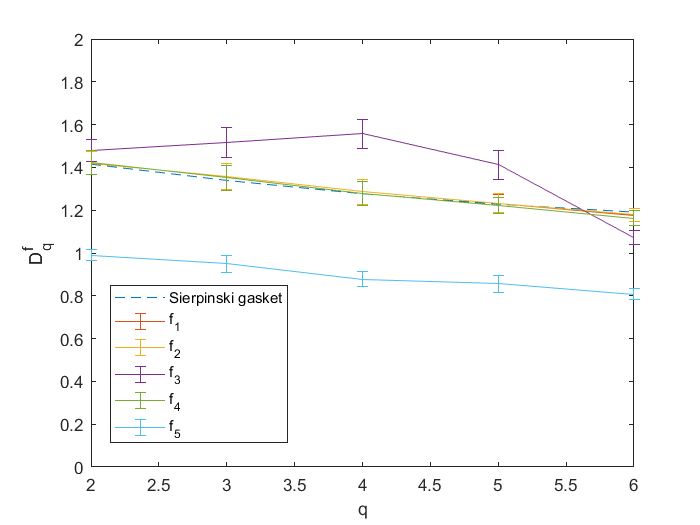}
\caption{Numerical estimates of $D_q^f$ for different observations: $f_1=Id$, $f_2(x,y)=(2x+y,2y)$, $f_3(x,y)=(\sin(\frac1x),\cos(\frac1y))$, $f_4(x,y)=((x-0.5)^2,2y)$ and $f_5=(1,y^2+x)$. In dashed lines is the $D_q$ spectrum of the underlying system. Estimates are computed as described in the text.}
\label{1}
\end{figure}

In \cite{obsrec}, we show that for the baker's map, a typical linear unidimensional projection gives $D_q^f=1$ for all $q$. Overall, this result along with our numerical computations suggest that Hunt and Kaloshin's results may extend to $q>2$ for a certain class of measures and for some smooth observations.
 
\section{The Extremal Index $\theta^f_q$}

When we work with real trajectories (i.e. when $f=Id$), the Extremal Index $\theta_q$, and more specifically the quantity $$H_q=\frac{\log(1-\theta_q)}{1-q},$$ encodes the hyperbolic properties of the system (see \cite{dq} for a detailed review). In particular, $H_q$ as a function of $q$ is constant for maps with constant Jacobian and is close to the metric entropy of the system (its Lyapunov exponent in dimension 1). When we introduce an observation $f$, the use of the Extremal Index to quantify the rate at which nearby trajectories diverge becomes less relevant, partly due to the fact that two nearby points in observational space may have antecedents far away in the actual phase space of the system. Let us investigate more in details this matter.\\

Keller and Liverani \cite{kl} provide a general formula for the Extremal Index of time series originating from dynamical systems. Applied to the present situation, and if the limits defining the different quantities exist, we have that

\begin{equation}\label{deftheta}
\theta^f_q=1-\sum_{k=0}^{\infty}p_{k,q},
\end{equation}

where 

\begin{equation}\label{p0}
p_{0,q}=\lim_{n\to\infty}\frac{\mu_q(S^q_n \cap T^{-1} S^q_n)}{\mu_q(S^q_n)}
\end{equation}

and for $k\ge 1$,

\begin{equation}\label{pk}
p_{k,q}=\lim_{n\to\infty}\frac{\mu_q(S^q_n \cap \bigcap_{i=1}^k T^{-i}(S^q_n)^c \cap T^{-k-1} S^q_n)}{\mu_q(S^q_n)}.
\end{equation}

In this general set up, obtaining a formula for $\theta^f_q$ is challenging, so let us place ourselves in the more simple case of expanding maps of the unit interval $I=[0,1]$. We define the following sets for $x\in I$ :
$$A_0(x)=\{y \in I, f(y)=f(x) \cap f(Ty)=f(Tx)\}$$
and
$$A_k(x)=\{y \in I, f(y)=f(x), f(T^iy)\neq f(T^ix),\text{ for } i=1,..,k\text{ and } f(T^{k+1}(y))=f(T^{k+1}(x))\}.$$

\begin{proposition}
Let $T$ be an expanding map of the unit interval $I=[0,1]$ which is $C^1$ by part and admitting an a.c. invariant measure $\mu(x)=h(x)dx$. Let $f : I \to J \subset \mathbb{R}$ be $C^1$ by part, finite to one and such that $f' \neq 0$ on $I$, then if

\begin{equation}\label{h1}
\mu(\{x \in I,A_0(x)=\{x\} \})=1
\end{equation} 

and, for all $k\ge 1$,

\begin{equation}\label{h2}
\mu(\{x\in I, A_k(x)= \emptyset \})=1,
\end{equation}

 we have that
 
\begin{equation}\label{thetaq}
\theta^f_q=1-\frac{\int_I \frac{h(x)^q}{\max(|f'(x)|,|(f\circ T)'(x)|)^{q-1}}dx}{\int_I \sum_{(y_1,...y_{q-1})\in (f^{-1}\{f(x)\})^{q-1}} \prod_{i=1}^{q-1}\frac{h(y_i)}{|f'(y_i)|}h(x)dx}.
\end{equation}
\end{proposition}

For a given map $T$, assumptions \ref{h1} and \ref{h2} should be satisfied for a generic observation $f$. The cases where these assumptions are not satisfied are when $T$ and $f$ share some particular symmetries and similarities of structures. For example, $\mu(A_0(x) = \{x\})\neq 1$ if both the graphs of $T$ and $f$ are symmetric with respect to the straight line of equation $x=1/2$.\\

\begin{proof}
We will write it for $q=2$. Following the lines of the proof in \cite{synchro} (where the case $f=Id$ is treated), and making use of the mean value theorem, we get:
\begin{equation}\label{1}
\begin{aligned}
\mu_2(S^2_n)&\sim \int_I \sum_{y\in f^{-1}\{f(x)\}} \mu(B(y,\frac{e^{-u_n}}{|f'(y)|})) d\mu(x)\\
                     &\sim 2e^{-u_n} \int_I \sum_{y\in f^{-1}\{f(x)\}} \frac{h(y)}{|f'(y)|} h(x)dx.
\end{aligned}
\end{equation}

We also have

\begin{equation}\label{2}
\begin{aligned}
\mu_2(S^2_n \cap T^{-1} S^2_n) &\sim \int_I \sum_{y\in A_0(x)} \mu(\{z\in I, z\in B(y,\frac{e^{-u_n}}{|f'(y)|})\cap Tz \in B(Ty,\frac{e^{-u_n}}{|f'(Ty)|}\})d\mu(x)\\
                                                   &\sim \int_I \sum_{y\in A_0(x)} \mu(\{z\in I, |z-y| \le \frac{e^{-u_n}}{|f'(y)|} \cap T'(y)|y-z| \le \frac{e^{-u_n}}{|f'(Ty)|}\}) h(x)dx.\\
                                                   &= \int_I \sum_{y\in A_0(x)} \mu(\{z\in I, |z-y| \le  \min(\frac{e^{-u_n}}{|f'(y)|},\frac{e^{-u_n}}{|T'(y)f'(Ty)|})\}) h(x)dx.\\
                                                   &\sim  2e^{-u_n} \int_I \sum_{y\in A_0(x)} \frac{h(y)h(x)}{\max(|f'(y)|,|(f\circ T)'(y)|)} dx.\\
\end{aligned}
\end{equation}

By a similar reasoning, we get that for $k \ge 1$,

\begin{equation}\label{3}
\mu_2(S^2_n \cap \bigcap_{i=1}^k T^{-i}(S^2_n)^c \cap T^{-k-1} S^2_n) \sim 2e^{-u_n} \int_I \sum_{y\in A_k(x)} \frac{h(x)h(y)}{\max(|f'(x)|,|(f(T^{k+1}(y))'|)}dx.\\
\end{equation}

Finally, combining eqs. \ref{deftheta},\ref{1}, \ref{2} and \ref{3}, we obtain

\begin{equation}
\theta^f_2=1 - \sum_{k=0}^{+\infty} \frac{\int_I \sum_{y\in A_k(x)} \frac{h(x)h(y)}{\max(|f'(x)|,|(f\circ T^{k+1})'(y)|)}dx}{\int_I \sum_{y\in f^{-1}\{f(x)\}} \frac{h'(y)}{|f'(y)|} h(x)dx}.
\end{equation}

This formula is still difficult to handle, but under condition \ref{h2}, we have that $p_{k,2}=0$ for $k>0$, and if moreover condition \ref{h1} holds, we obtain 

\begin{equation}\label{thetafin}
\begin{aligned}
\theta^f_2&=1-p_{0,2}\\
             &=1-\frac{\int_I \frac{h(x)^2}{\max(|f'(x)|,|(f\circ T)'(x)|)} dx}{\int_I \sum_{y\in f^{-1}\{f(x)\}} \frac{h(y)h(x)}{|f'(y)|}dx}.
\end{aligned}
\end{equation}

We can generalize this result for $q\ge 2$ to obtain the desired result.
\end{proof}

\begin{example} Let us take $Tx=2x \mod 1 $ and $$f(x)=\left\{
    \begin{array}{ll}
        2x & \mbox{ if } 0\le x\le 1/2 \\
        3/2-x & \mbox{ if } 1/2<x\le 1.\\
    \end{array}
\right.$$

The function $f$ is not injective, satisfies conditions \ref{h1} and \ref{h2}, and computations are worked out quite easily, which constitutes a good test for our results. Applying formula \ref{thetaq}, we get $$\theta^f_q=1-p_{0,q}=1-\frac{2+2^{2-q}}{1+3^q}.$$ This result is confirmed by numerical experiments (see figure \ref{1}). We used the estimate $\hat{\theta}_{5}$ introduced in \cite{ei}. This estimate consists in evaluating by means of Birkhoff sums the 5 first $p_{k,q}$ terms introduced in formula \ref{pk} and subtracting them from 1. It requires to fix a high threshold $u$ (here we take $u$ equal to the 0.99999-quantile of the $Y_i$ distribution). As expected, we find that all the $p_{k,q}$ are 0 or very close to 0 for $k\ge1$. The result is displayed in figure \ref{1}. Our results are averaged over 10 runs, with trajectories of length $2.10^7$. We see easily that in this example, $D_q^f=1$ for all $q$.\\

\end{example}

\begin{figure}[h!]
\centering
\includegraphics[height=3.5in]{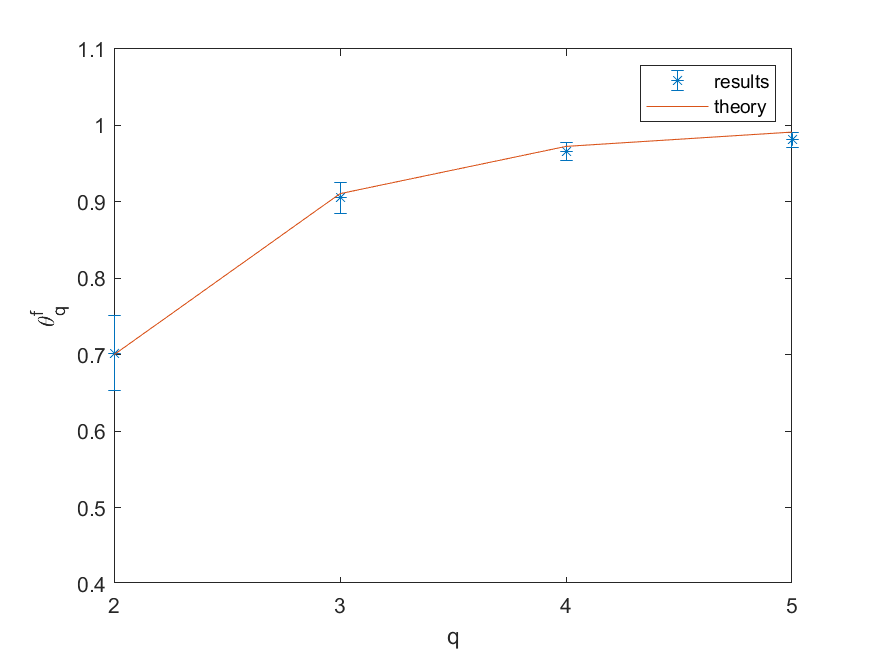}
\caption{Comparison between theory and computation for the $\theta_q^f$ spectrum of the system described in the text.}
\label{1}
\end{figure}

A general formula for higher dimensional system is out of scope, but we expect that with conditions of `non compatibility' between the dynamics and the observation analogue to conditions \ref{h1} and \ref{h2}, all the $p_{k,q}$ terms are 0 for $k\ge1$. This hypothesis is verified by several numerical experiments, as we will see. The vanishing of the $p_k$ terms is particularly welcome when it comes to proving the existence of the extreme value law using more classical approaches in Extreme Value Theory \cite{book,freitascond1,freitascond2}.\\

The presence of the derivative of the observation in formula \ref{thetaq} renders the interpretation of $\theta_q^f$ more difficult than in the case $f=Id$, however we notice two facts :

\begin{itemize}
\item For a given observation $f$, the larger are the values of $|T'|$ over phase space, the larger are the values of $\theta_q^f$, so this index can somehow still quantify the hyperbolic properties of $T$.
\item For a given map $T$, the more the points in the observational space have antecedents by $f$, the bigger is the denominator in equation \ref{thetaq}, and the larger is $\theta_q^f$: oscillatory observations give higher values for the extremal index.
\end{itemize}

We expect analogous properties to hold for higher dimensional systems. We now test that statement.\\

 In figure \ref{thetf} (a), we compare the estimates of $\theta_q^f$ for the 2-dimensional H\'enon system, defined by $T(x,y)=(1-ax^2+y,bx)$, with $a=1.4$ and different values of $b$ such that the system admits a strange attractor. The observation we take is $f(x,y)=\frac{x+y}{2}$. The determinant of the Jacobian is given by $b$. We find indeed that for this fixed choice of observation, the more the original system tends to separate trajectories (the higher is parameter $b$), the higher are the values of $\theta_q^f$, even for lower dimensional projections. The estimates $\hat{p}_{k,q}$ of the $p_{k,q}$ terms, for $k>0$ are all null or close to 0 for all the observation that we considered, as conjectured before.\\

In figure \ref{thetf} (b), we plot the estimates of the extremal index for 2-dimensional H\'enon system (with usual parameters $a=1.4$, $b=0.3$) and different observations. We observe that for one to one observations, ($f_1$, $f_2$ and $f_3$), the $\theta_q^f$ spectrum remains low, although the form of the Jacobian can impact significantly the values of $\theta_q^f$. When the observation ceases to be one to one, the spectrum of extremal indices increases significantly (see the curve for $f_4$). This effect is even more important for the very oscillatory function $f_5$.

%

\begin{figure}
\centering
\begin{subfigure}{.5\textwidth} \label{thet2}
  \centering
  \includegraphics[width=1\linewidth]{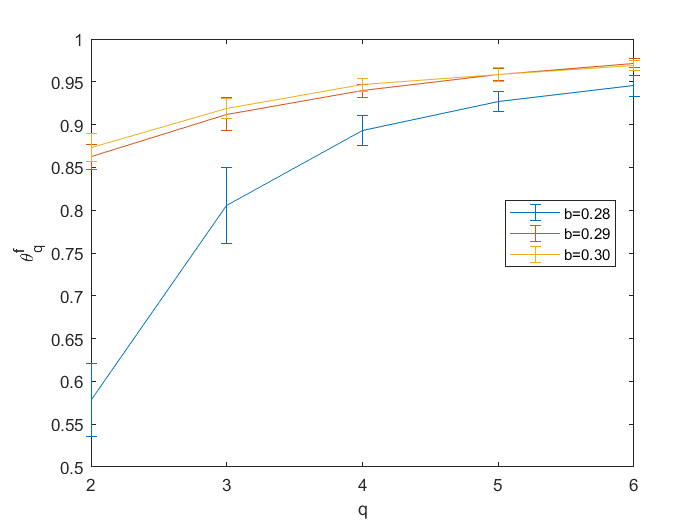}
  \caption{ }
\end{subfigure}%
\begin{subfigure}{.5\textwidth} \label{thetf}
  \centering
  \includegraphics[width=1\linewidth]{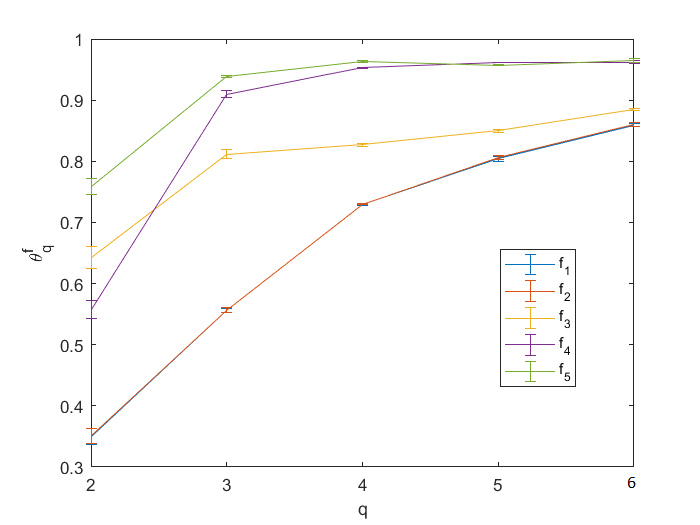}
  \caption{ }
\end{subfigure}
\caption{Left: Estimates for the $\theta^f_q$ spectrum computed for a H\'enon system with different parameters $b$ and for the observation $f(x,y)=\frac{x+y}{2}$. Right: Estimates for the $\theta^f_q$ spectrum computed for the H\'enon system (b=0.3) and different observations : $f_1=Id$, $f_2(x,y)=(100x+y,100y)$, $f_3(x,y)=(x,100y)$, $f_4(x,y)=(x^2,y^2)$, $f_5(x,y)=(\sin(1/x),\cos(1/y))$. For both figures, we used the estimate $\hat{\theta}_5$, with trajectories of length $10^6$ and a threshold value equal to the $0.999$ quantile of the $\phi^f_q$ distribution. The error bars represent the standard deviation of the results over 10 runs.}
\label{thetf}
\end{figure}

%

\section{Acknowledgement}
The author would like to thank J\'er\^ome Rousseau and Sandro Vaienti for fruitful discussions concerning this work.

\end{document}